\documentclass[12pt,reqno]{amsart}
\usepackage{amssymb,latexsym,amstext,amsfonts,amscd,amsmath, }
\usepackage{graphics}
\usepackage{epsfig}
\theoremstyle{plain}
\newtheorem{theorem}{Theorem}[section]

\newtheorem{lemma}[theorem]{Lemma}

\newtheorem{proposition}[theorem]{Proposition}
\newtheorem{corollary}[theorem]{Corollary}

\numberwithin{equation}{section}

\theoremstyle{remark}
\newtheorem{rema}[theorem]{Remark}
\newtheorem{exam}[theorem]{Example}

\renewcommand{\labelenumi}{(\roman{enumi})}
\makeatletter
\def\alphenumi{%
        \def\theenumi{\alph{enumi}}%
        \def\p@enumi{\theenumi}%
        \def\labelenumi{(\@alph\c@enumi)}}
\makeatother

\newcommand{\CC}{\mathbb{C}}
\newcommand{\C}{\mathbb{C}}
\newcommand{\RR}{\mathbb{R}}
\newcommand{\R}{\mathbb{R}}

\newcommand{\NN}{\mathbb{N}}

\newcommand{\I}{\mathbb{I}}
\newcommand{\GG}{\bold G}

\newcommand{\p}{\partial}
\newcommand{\dsp}{\displaystyle}
\begin{document}
\pagestyle{plain}

\title{ Hyperbolic polynomials and multiparameter real analytic perturbation
theory}
\author{Krzysztof Kurdyka}
\address{
Laboratoire de Mathematiques (LAMA), Universit\'e de Savoie\\UMR 5127
CNRS\\ \newline  73-376 Le Bourget-du-Lac cedex FRANCE }
\email{ Krzysztof.Kurdyka@univ-savoie.fr}
\author{Laurentiu Paunescu}
\address{School of Mathematics and Statistics\\
University of Sydney, NSW 2006, Australia}
\email{laurent@maths.usyd.edu.au}
%\urladdr{}
\thanks{The first author thanks University of Sydney for
support. Part of this work was done while the second author was
visiting MSRI Berkeley.}
\dedicatory{}
\date{21 January 2006}
%\maketitle
\begin{abstract}

Let
$P(x,z)= z^d +\sum_{i=1}^{d}a_i(x)z^{d-i}$ be a polynomial, where  $a_i$
are real
analytic functions  in an open subset $U$ of $\R^n$.
If for any
$x \in U$
the polynomial
$z\mapsto P(x,z)$ has only real roots, then we can write those roots as
locally lipschitz functions of $x$. Moreover, there exists a modification
    (a locally finite composition of blowing-ups with smooth centers)
    $\sigma : W \to U$
    such that the roots of the  corresponding polynomial
    $\tilde P(w,z)  =P(\sigma (w),z),\,w\in  W $, can be written
locally as analytic functions of $w$.
Let $A(x), \, x\in U$ be an analytic family of symmetric matrices, where
$U$  is  open in
$\R^n$. Then there exists a modification
    $\sigma : W \to U$, such the corresponding family
$\tilde A(w) =A(\sigma(w))$ can be locally diagonalized analytically
(i.e. we can choose locally a  basis of eigenvectors in an analytic way).
This generalizes the Rellich's  well known  theorem (1937) for one parameter
families. Similarly  for an analytic family $A(x), \, x\in U$ of
antisymmetric matrices there exits a modification $\sigma$ such that we
can find locally  a basis of proper subspaces in an analytic way.
\end{abstract}

\maketitle
\section{Introduction}

In the late 30's F. Rellich \cite{rellich},\cite{rellichbook} developed
the theory of one parameter analytic perturbation theory of linear operators.
This theory culminates  with the celebrated monograph of T. Kato \cite{kato}.
To study the behaviour of eigenvalues of symmetric matrices under analytic
one-parameter perturbation Rellich  proved the following fundamental
fact.
Let \begin{equation}\label{maineq}
P(x,z)= z^d +\sum_{i=1}^{d}a_i(x)z^{d-i}
\end{equation} be a
polynomial, where
$a_i$ are real
analytic functions on an open interval $I\subset \R$. If
for any
$x \in I$
the polynomial
$z\mapsto P(x,z)$ has only real roots (we call such a polynomial
{\it hyperbolic}), then there are analytic functions
$f_i:I\to \R, \,i=1,\dots ,n$ such that
$P(x,z) =\prod_{i=1}^d (z-f_i (x))$. In other words we can choose
analytically  the roots of $P$.

If we consider a multiparameter version of this theorem, i.e. we
assume now that  $a_i$
are real
analytic functions in an open subset $U$ of $\R^n$, $n>1$,
    then
we have  a simple counterexample $z^2 -(x_1^2 +x_2^2)$.
For this reason a multiparameter perturbation theory was not developed (to
our knowledge), though it was  suggested  by Rellich. In this paper we
give some  generalizations of Rellich's theory in the multiparameter
case. These generalizations are purely real,  they  make no sense in
the complex case developed by Kato.

The first generalization was inspired by S. {\L}ojasiewicz, who suggested
that the roots of the polynomial $P$ can be chosen  locally in a
lipschitz way. This is true as we prove in Theorem \ref{liphip}. The
result  is quite delicate since  in the one
parameter case (the Rellich Theorem) there  is no way  to bound  the 
lipschitz constant
for roots in terms of bounds for the coefficients $a_i(x)$. The proof 
is obtained
by a reduction to  the $2$-parameter case and  by a careful study of
a desingularization of singularities of the zeros of
$P$. In fact we are able  to keep track of partial derivatives of roots after
blowing up, because we are dealing with a family of hyperbolic
polynomials. This is rather surprising, since it is known that this is
impossible in  general, for instance there are
blow-analytic (or arc-analytic functions)  which are not locally lipschitz.
Our result is related to  Lidski's Theorem which implies that
the spectral mapping on the space of symmetric matrices is lipschitz
(globally). More precisely,  our Theorem \ref{liphip} implies the above
corollary of Lidskii's theorem (in a weaker form), however  Lidskii's theorem
does not imply our theorem. Indeed as Rellich noticed in
\cite{rellichbook} not every analytic family of hyperbolic
polynomials can be written as  characteristic polynomials of an analytic
family of symmetric matrices. Surprisingly this is related to Hilbert's 17th
   problem for analytic functions. In a similar way we prove an
    analogue  of  Theorem
\ref{liphip} in the case where  all the roots  of   $z\mapsto P(x,z)$ are
purely imaginary (we call such a polynomial {\it antihyperbolic}). This
is important in the study of analytic families of antisymmetric
matrices.

The second direction of generalization is related to the theory of
arc-analytic functions; initiated by the first author in \cite{kus}.
Actually the roots of an analytic family of hyperbolic polynomials can be seen
as a multivalued arc-analytic function.  As we prove in Theorem
\ref{hipres} it  turns out that after  suitable blowing-ups of the space
of parameters we can write locally  the roots  of hyperbolic polynomials as
analytic functions of parameters.

Rellich's theory deals not only with eigenvalues but also with
eigenvectors. He proved (see for instance \cite{rellichbook}) that every
one parameter analytic  family of symmetric matrices admits  an analytic
choice of bases of eigenvectors. In other words such a family can be
analytically diagonalized even if the eigenvalues become multiple. As we
prove in Theorem \ref{diagres} this can be also done for   multiparameter
analytic families, but first we have to blow up the space of parameters in
order to make the eigenvalues  locally normal crossing.

Finally, we also study  analytic families of antisymmetric matrices
depending on several parameters. We prove analogously that, after
suitable blowing-ups
of the parameters, we can
reduce them locally    to the canonical form  in an analytic way.

\section{ Arc-analytic functions}\label{arcan}
For further convenience we recall here some facts concerning arc-analyticity.
Let $U$ be an open subset of  $\RR^n$.  Following \cite{kus} we say
that a map
$f:U\to\RR^k$ is {\it arc-analytic} if for any analytic arc
$\alpha:(-\varepsilon, \varepsilon) \to U$, the composed function
$f\circ\alpha$ is also analytic.

In general arc-analytic maps are very far from being analytic, in particular
there
are arc-analytic functions  which are not subanalytic \cite{ku1}, not
continuous
\cite{bmp},  with a non-discrete singular set \cite{ku2}. Hence it is
natural
to consider
only arc-analytic maps with subanalytic graphs.
Earlier T.-C. Kuo \cite{Kuo85}, motivated by  equisingularity problems,
introduced     the notion
of {\it blow-analytic} functions, i.e. functions which become
analytic  after
a composition with
appropriate  proper bimeromorphic maps (e.g. a  composition of blowing-ups
with
smooth centers).
Clearly any blow-analytic mapping is arc-analytic and subanalytic.
The converse holds in a slightly weaker form
\cite{Bierstone-Milman90}, see also \cite{Par}. We shall explain it in the next
section.

Blow-analytic maps have been
studied by several authors (see the survey \cite{kkf}).
It is known
that  in general subanalytic and  arc-analytic functions are continuous
\cite{kus}, but  not necessarily (locally)  lipschitz  \cite{kkf},
\cite{Paunescu}.

The following examples are arc-analytic but not analytic functions:
$$
f =\frac{x^3}{x^2+y^2} , \, g = \frac{xy^5}{x^4+y^6},
\, h = \sqrt{x^4+y^4}.
$$

The function $f$  is locally lipschitz, but not $C^1$ (cf. \cite{kus}),
the function $g$  is not locally lipschitz (cf. \cite{Paunescu}). The
function $h$   is $C^1$ but not $C^2$ (cf. \cite{Bierstone-Milman90}).

Recently  we have proved in \cite{kurpalip} that, if $h$ is 
arc-analytic and $h^r$ is
analytic for some integer
$r$, then $h$ is locally lipschitz. However arc-analytic roots  of
polynomials with analytic coefficients are not necessarily lipschitz.
\begin{exam}\label{exem1}\rm Consider a polynomial
$P(x,y,z) = (z^4 -(x^2+y^8))^2 -x^4-y^{20}$. It has an arc-analytic
%(cf. Section \ref{arcan})
root
$$
f=\root 4\of{{x^2+y^8} -\sqrt{x^4+y^{20}}},
$$

which is not lipschitz! Note that the above polynomial is not hyperbolic.
\end{exam}
It is useful to consider  arc-analytic  complex
valued functions, where we understand that they are analytic on real
analytic arcs.
We cannot avoid arc-analytic solutions in the sense above. Indeed we
have the following type of examples which appear in our context:

\begin{exam} Consider
$P(z,x,y)=z^4-x^8-y^8$ as polynomial in $z$, then it has   the obvious
roots:

$$ z_1= \root 4\of{x^8+y^8}, z_2=- \root 4\of{x^8+y^8},
z_3=i \root 4\of{x^8+y^8}, z_4=-i \root 4\of{x^8+y^8}.$$
\end{exam}

\subsection{Locally blow-analytic functions}

We recall some of the notions used in this paper
(for more information see for instance \cite{bmp}, \cite{kkf}, \cite{fkp},
\cite{Kuo85}, \cite{ku1}, \cite{ku2}, \cite{pl}).

We recall first a definition of a local blowing-up. Let $M$ be  an analytic
manifold and
$\Omega
\subset M$  an open set. Assume that $X$ is an analytic  submanifold of $M$,
closed in
$\Omega$. Then we can define
   $\tau:\tilde \Omega \to \Omega$,  the
blowing-up of
$\Omega$ with the center  $X$, see for instance \cite{H73}
or \cite{loj}.
A restriction of
$\tau$ to an open subset of $\tilde \Omega$ is called {\it a local
blowing-up with
a smooth (nowhere dense) center}.

Let $U$ be a neighbourhood of the origin of $\RR^n$ and
let $f:U \to\RR^m$ denote a map defined on $U$ except
possibly some nowhere dense subanalytic subset of $U$.
We say that $f$ is {\it locally blow-analytic} via a locally finite
collection of analytic modifications $ \sigma_{\alpha}: W_{\alpha} \to
\RR^n$, if for each $\alpha$ we have
\begin{enumerate}
\item $W_{\alpha}$ is isomorphic to $\RR^n$ and $\sigma_{\alpha}$ is the
composition
of finitely many local blowing-ups
with smooth nowhere dense centers, and $f \circ \sigma_{\alpha}$ has an
analytic extension on $ W_{\alpha}$.
\item There are subanalytic compact subsets  $K_{\alpha} \subset  W_{\alpha}$
such that $\bigcup \sigma_{\alpha}(K_{\alpha})$ is a neighbourhood  of
$\overline U$.
\end{enumerate}
The notion of  {\it (locally) blow-analytic} functions (or maps) is very
much related to the notion of {\it arc-analytic} functions.
Indeed in \cite{Bierstone-Milman90}, see also \cite{Par}, it is proved that
an {\it arc-analytic} function has {\it  subanalytic} graph if and
only if it is  {\it locally blow-analytic}.
\begin{rema} .
The definition of arc-analytic function is  much more intrinsic and it is
usually easier to
check that a given function is  arc-analytic  than to check that it is
blow-analytic.
Actually, when the first author introduced (in mid 80's) arc-analytic
functions he has
hoped that subanalytic and arc-analytic are
exactly the  same with (globally)  blow-analytic.
This is true for semialgebraic functions and for functions in $2$ variables
(since we blow up only points).
In a forthcoming paper (joint with A. Parusi\'nski) the authors
shall give a
proof of this conjecture for
functions in $3$ variables. But the general case presents serious
difficulties and remains still open.
\end{rema}

\medskip

\section{Hyperbolic Polynomials}\label{ssplitting}
\medskip\medskip

\subsection{Splitting lemma for polynomials}

Given a $p$-tuple $a=(a_1,\dots , a_p) \in \R^p$ and
     a $q$-tuple $b=(b_1,\dots , b_q) \in \R^q$, we associate
two polynomials
$$
P_a(z)= z^p +\sum_{i=1}^{p}a_iz^{p-i}, \,\,
%\text{ and }
Q_b(z)= z^q +\sum_{j=1}^{q}b_jz^{q-j}.
$$
We consider  the product of these polynomials
     $$\displaystyle P_aQ_b =R_c= z^{p+q} +\sum_{k=1}^{p+q}c_kz^{p+q-k},
$$ where
$c=(c_1,\dots , c_{p+q}) \in \R^{p+q}$. This defines a polynomial map
$$\Phi: \R^p\times \R^q \ni (a,b) \mapsto c \in \R^{p+q}.$$
     The following lemma is crucial.

\begin{lemma}{\bf (Hensel's Splitting Lemma)}\label{splitting}
     \begin{enumerate}
\item The jacobian of $\Phi$ at $(a,b)$ is equal (up to sign) to the
resultant of $P_a$ and $Q_b$.
\item Let us fix $\bar a=(\bar a_1,\dots , \bar a_k) \in \R^p$ and
$\bar b=( b_1,\dots ,\bar b_q) \in
\R^q$   and assume that corresponding polynomials $P_{\bar a}$ and
$Q_{\bar b}$ have no common zeros
in  $\CC$, in other words that their resultant is non zero.
Then there exists  a neighbourhood $U \subset \R^{p+q}$ of $\bar c =
\Phi(\bar a, \bar b)$ such that
     for any $c \in U$ the corresponding polynomial  splits in a
unique way,  $R_c =  P_a Q_b$, moreover
the mapping $a= a(c), b= b(c)$ is analytic (even Nash) and satisfies
$\bar a= a(\bar c), \bar b= b(\bar
c)$ .
\end{enumerate}
\end{lemma}

Indeed, it is easy to see that the Jacobian matrix of $\Phi$  is exactly
the Sylvester matrix of the
pair
$P_a$ and $Q_b$. So the resultant  $Res(P_a,Q_b)$,  which is by the
definition   equal to (up
to sign) the determinant  of the Sylvester matrix,   is also equal to  the
jacobian of
$\Phi$ at
$(a,b)$, see e.g. \cite{Bierstone-Milman90},\cite{a}. Recall that two
polynomials $P_{ a}$ and
$Q_{ b}$ have no common zeros in  $\CC$, if and only if  their resultant is
non zero. The
second  part of the lemma is just a consequence of the fact that
$\Phi$ is invertible in a neighbourhood of $(a,b)$, in particular $\Phi^{-1}$
is analytic (even Nash) by
the Inverse Mapping Theorem.

In the sequel we will use the following consequence of the splitting Lemma.
\begin{corollary}\label{analsplitting}
Let $R(x,z)= z^r +\sum_{k=1}^{r}c_k(x)z^{r-k}$, where $c_k(x)$ are analytic
functions
in some open set   $ \Omega \subset \R^m$. Assume that for some $x_0\in
\Omega$ the polynomial
$z\mapsto R(x_0,z) $ splits,  i.e.  $R(x_0,z)= P_{x_0}(z) Q_{x_0}(z) $,
where $\deg P_{x_0} =p$,
$\deg Q_{x_0} =q$   and $r=p+q$. Suppose moreover that $P_{x_0}(z)$ and
$ Q_{x_0}(z) $ have no
common roots in $\CC$. Then there exist a neighbourhood $U \subset \Omega$
of $x_0$, and analytic
functions,  $a_i:U\to \R,\, i=1,\dots, p$ and $b_j:U\to \R,\, j=1,\dots, q$
such that
$$
R(x,z) =P(x,z) Q(x,z), \, x\in U, \, z\in \R,
$$
where $P(x,z)= z^p +\sum_{i=1}^{p}a_i(x)z^{p-i}$, $Q (x,z)= z^q +
\sum_{j=1}^{q}b_j(x)z^{q-j}$.
Moreover $P(x_0,z) =P_{x_0}(z)$ and  $Q(x_0,z) =Q_{x_0}(z)$.

\end{corollary}
{\it Remark.} Splitting Lemma  and Corollary \ref{analsplitting} hold of
course over complex numbers,
but we don't need this.

\subsection{Newton-Puiseux Expansions}

For latter use we  recall some  classical facts  about the roots of Weierstrass
polynomials.
Let
\begin{equation}\label{classic}P(x,z)= z^d +\sum_{i=1}^{d}a_i(x)z^{d-i},
\end{equation}
with $a_i$
real analytic functions in a neighbourhood of $0\in \R$.
%, be a  Weierstrass polynomial.
Then there are  holomorphic functions $h_i, \, i=1,\dots, d$ and an integer
$r$ such that
$$P(x,z)=\prod_{i=1}^d [z-h_i (x^{1/ r})]=\prod_{i=1}^d [z-f_i (x)]$$
for $x\ge0$ close enough to $0$, and any $z\in\CC$.
We call $f_i(x)=h_i (x^{1/ r})$ a {\it Newton-Puiseux root of $P$}. Clearly
each $f_i$ is given by a {\it Puiseux expansion} $f_i(x) =
\sum_{\nu=0}^\infty \alpha^i_\nu x^{\nu/r}$.

\subsection{Hyperbolic polynomials}

Let $$P(z)= z^d +\sum_{i=1}^{d}a_iz^{d-i}$$ be a polynomial with real
coefficients.
Let $z_1,\dots, z_d$
be all complex roots    of $P$;  recall that
$a_1=z_1+\dots+ z_d$.
By {\it Tschirnhausen transformation}, which is the change of variable
$z\mapsto z -\frac{a_1}{d}$, we may assume that
$a_1= 0$.
We say that
$P$ {\it is hyperbolic} if all its roots are real. Hyperbolic polynomials
appear naturally,
     for instance as  characteristic polynomials of symmetric matrices.
We state now two elementary but crucial properties of hyperbolic polynomials.
\begin{lemma}\label{hyproots}
Let  $P(z)= z^d +\sum_{i=2}^{d}a_iz^{d-i}$ be a polynomial with real
coefficients (note that
$a_1=0$).  Denote  the roots (possibly  complex) of $P$ by
$z_1,\dots, z_d$ . Then
\begin{equation}\label{hyprootseq}
z_1^2+\dots+ z_d^2 = -2a_2.
\end{equation}
     Consequently, if $P$ is hyperbolic, then $a_1=a_2  =0$ if and 
only if  $P(z) =z^d$,
that is $0$ is  the only
root of $P$.
\end{lemma}

Proof:  since $a_2=\sum_{i<j}z_iz_j$, we have $z_1^2+\dots+ z_d^2 =
a_1^2 -2a_2=-2a_2$.
If all $z_i$ are real, then $z_1^2+\dots+ z_d^2=0$ implies that all $z_i=0$.

In the sequel we will study
families of monic polynomials depending analytically on parameters
(i.e. the coefficients are analytic functions of
parameters), such that for each values of the parameters the
corresponding polynomial
is hyperbolic. We will also
call,  for short, such a family {\it a hyperbolic polynomial}.

\subsection{Rellich's Theorem}
In the late  30's Rellich \cite{rellich} proved a rather surprising
fact about the
roots of hyperbolic
polynomials of the form
\eqref{classic}. His result is the following.
\begin{theorem}\label{rellich}{\bf(Rellich 1937)}
Consider a polynomial
\begin{equation*}P(x,z)= z^d +\sum_{i=1}^{d}a_i(x)z^{d-i},
\end{equation*} with $a_i$
real analytic functions in an open interval $I\subset \R$. Assume that for
each $x\in I$ all the roots of the polynomial $z\mapsto P(x,z)$ are real. Then
there exist real analytic
functions $f_i:I\to \R$ such that
\begin{equation}\label{rellicheq}
P(x,z)= \prod_{i=1}^d [z-f_i (x)], \,x\in I, z\in \R.
\end{equation}

\end{theorem}
We will outline   a proof of the above theorem.

It is inspired by \cite{a}, however we made it
shorter since  we use Puiseux's theorem.

Note that by the analytic extension argument it is enough to prove the theorem
locally. We fix a
point, say $0\in I$, and  assume that $a_i$ are analytic in a neighbourhood of
$0$.

\begin{itemize}
\item
     1st Step:

We may assume that all $a_i(0) =0$.

Indeed, if $P(0,z)= (z-c)^d$ then, shifting $z\mapsto z-c$, we may assume
that $c=0$.
Consequently all $a_i(0) =0$.
Otherwise $P(0,z) =(z-c)^pP_2(z)$, $0<p<d$, with $P_2(c)\ne 0$.
     Applying Corollary \ref{analsplitting}
we can split our polynomial as $P(x,z)= P_1(x,z) P_2(x,z)$, where $P_1$ and
$P_2$ are of the form
\eqref{classic} with real analytic coefficients in a neighbourhood of
$0\in \R$.
Hence we can handle separately $P_1$ which is already of the form considered
above, and  $P_2$
which is of a smaller degree.
\item
2nd Step:

     Let us write $$P(x,z)= z^d +\sum_{i=1}^{d}a_i(x)z^{d-i} =\prod_{i=1}^d
[z-f_i (x)], \, x>0,$$
where all $a_i$ are real analytic in a neighbourhood of $0\in \R$, and
$a_i(0) =0$, hence also
$f_i(0) = 0$. Applying Tschirnhausen transformation $z\mapsto z -
\frac{a_1(x)}{d}$ we may assume that
$a_1(x)\equiv 0$.
Denote by  $z_1(x),\dots, z_d(x)$ all the roots of $z\mapsto P(x,z)$. Then
Lemma \ref{hyproots} yields
\begin{equation}\label{newton}
-2a_2(x) = z_1(x)^2+\dots+ z_d(x)^2,
\end{equation}
for $x$ in a neighbourhood of $0\in \R$.
But by our assumption all the roots $z_i(x)$ are real, hence $a_2$ must be
negative.
Consequently the order of $a_2$ at $0$ is even and we can write
\begin{equation}\label{znak}a_2(x)= x^2b(x),
\end{equation}
with $b(x)$ analytic in a neighbourhood of $0\in \R$.
Applying \eqref{newton} and \eqref{znak} to the Puiseux roots of $P$ we obtain
\begin{equation}\label{klucz} f_1(x)^2+\dots+ f_d(x)^2= x^2b(x).
\end{equation}
So we easily deduce the following lemma.
\begin{lemma}\label{order1} The order of each $f_i$ at $0$ is greater or
equal than $1$.
\end{lemma}
By the order of $f_i$ we mean the smallest rational exponent in its Puiseux
expansion  such
that its coefficient  does not vanish.

Accordingly, by Vi\'ete's  formulas $\dsp a_i =(-1)^i\sum_{k_1< \dots <
k_i}z_{k_1}\dots  z_{k_i}  $, we obtain:

\begin{lemma}\label{orderi} The order of each $a_i$ at $0$ is greater or
equal than $i$.
\end{lemma}

\end{itemize}

Now we are in the position to conclude Rellich's theorem. We are going to show
that in the Puiseux expansion of each $f_i$ there are only integer exponents.\\
Let us   write
%contrary, more precisely that
$f_i(x) = \sum_{\nu=0}^\infty \alpha^i_\nu x^{\nu/r}$.
%and $a^i_k=\ne 0$ for some $k$ which is not a multiple of $r$, morover $k$
%is the smallest such an  integer.
By Lemma \ref{order1} we know that all
$ \alpha^i_1 =  \cdots =  \alpha^i_{r-1}=0$, so
$$
\frac{f_i(x)}{x} = \sum_{\nu=r}^\infty \alpha^i_\nu x^{\nu/r -1}
$$
are all bounded, and they are the Puiseux  roots of the polynomial
\begin{equation*}\tilde P(x,z)= z^d +\sum_{i=1}^{d}\frac{a_i(x)}{x^i}z^{d-i},
\end{equation*}
with $\tilde a_i(x)=\frac{a_i(x)}{x^i}$ real analytic at $0\in \R$, by
Lemma
\ref{orderi}. Now we apply the first step of the reduction to the polynomial
$\tilde P(x,z)$.  So  may assume that all $\tilde a_i(0)=0$. Note that
the shift affects only  the coefficient
$\alpha^i_r$. By Lemma
\ref{order1} we   deduce that $$\alpha^i_\nu =0, \, \nu=r+1,\dots,2r-1.$$

Continuing this  process we see that  for any integer $\nu$ which is not a
multiple of
$r$, the corresponding coefficient vanishes,  that is $\alpha^i_\nu =0$.
Actually  we can write $f_i$ as a convergent  power series
$f_i(x) = \sum_{n=0}^\infty \alpha^i_{rn} x^{n}$.

Formally our argument applies only for $x>0$, but since now we know that $f_i$
are analytic in a neighbourhood of $0\in \R$, we  deduce that also
for  negative $x$,   $f_1(x),\dots,f_d(x)$ are   the roots of the polynomial
$z\mapsto P(x,z)$.
So we can write
$$P(x,z)= z^d +\sum_{i=1}^{d}a_i(x)z^{d-i} =\prod_{i=1}^d [z-f_i (x)]$$
with $f_i$ analytic in a neighbourhood of $0\in\R$. By the analytic extension
argument, each $f_i$
extends to a unique analytic function on the whole interval $I$.
Hence Rellich's theorem follows.

\subsection{ Expansions  of the roots of  hyperbolic polynomials in 2
parameters}

     Consider a polynomial
\begin{equation}\label{manyvar}P(x,y,z)= z^d +\sum_{i=1}^{d}a_i(x,y)z^{d-i},
\end{equation}
with $a_i(x,y)$ analytic in a neighbourhood of $0\in \R^2$.

\begin{proposition}\label{keysplitprop}
     Assume that $P$ is hyperbolic with
respect to $z$, that is  for each $(x,y)$, the polynomial
$z \mapsto P(x,y,z)$ has only real roots. Then, in a  set
$H =\{\vert x\vert <y^N, \,0<y<\delta\}$, where $N$ is large enough
and $\delta $ is small enough, the polynomial $P$ splits in the form
\begin{equation}\label{keyspliteq}
P(x,y,z) =\prod_{i=1}^d [z-f_i (x,y)],
\, (x,y)\in H, \, z\in\R,
\end{equation}
with $f_i$ are analytic in $H$.
\end{proposition}
Now we can state a key proposition, which allows us to prove the fact that
the roots of hyperbolic polynomials are lipschitz.
\begin{proposition}\label{keylipprop}
Each function $y\mapsto \frac{\partial f_i}{\partial x}(0,y)$ extends to an
analytic function
in a neighbourhood of $0\in\R$, in particular $\frac{\partial f_i}{
\partial x}(0,y)$ are
bounded for $y\in (0,\delta)$.
\end{proposition}

{\it Proof of} Proposition \ref{keysplitprop}. We shall proceed by 
the  induction on the highest
multiplicity of  a root of the univariate polynomial $z \mapsto P(0,0,z)$.
\begin{itemize}
\item
  Case 0. All roots of $z \mapsto P(0,0,z)$ are simple, then the statement
of the proposition is an immediate consequence of the Implicit
Function Theorem.
\item
Case 1. Let  $c$ be a root of $z \mapsto P(0,0,z)$ of the maximal multiplicity,
We can suppose that $P(0,0,z)= (z-c)^d$.
Otherwise $P(0,0,z) =(z-c)^pP_2(z)$, $0<p<d$, with $P_2(c)\ne 0$.
     Applying Corollary \ref{analsplitting}
we can split our polynomial as $P(x,y,z)= P_1(x,y,z) P_2(x,y,z)$, where $P_1$
and $P_2$ are of the form
\eqref{manyvar}  with real analytic coefficients in a neighbourhood of $0\in
\R$.
Hence we can handle separately $P_1$, which is already of the form considered
above, and  $P_2$
which is of a smaller degree.
%%%%
Finally, shifting $z\mapsto z-c$, we may suppose  that
$c=0$.

     So in formula \eqref{manyvar}
we may suppose that
all $a_i(0,0)=0$.
\vskip 0,5 cm
Now we consider the hyperbolic polynomial $P(0,y,z)$. According to Rellich's
Theorem \ref{rellich}
we have
$$P(0,y,z) =\prod_{i=1}^d [z-c_0^i (y)]$$
with $c_0^i (y)$ analytic in a neighbourhood of $0\in \R$.
\item
Case 1.1.  Assume that not all $c_0^i (y)$ are identical as 
functions, note that
$c_0^i (0)=0$ for all $i$.  We are going to describe
an operation,
which will allow us to reduce the multiplicity of the root $c=0$.
Let $f,g$
be two  distinct analytic  functions in a neighbourhood of $0\in \R$,
then
$$f (y)-g (y)=y^k b(y)$$
where $b(y)$ is analytic and  $b(0)\ne 0$. We will call $k$ the {\it order of
contact of} $f$ and  $g$
at $0$.
%
%The  operation of {\it taking proper transform by blowing-up} decreases by 1
%the  order of  contact
%between distinct roots $c_0^i (y)$.
%So after taking proper transform finitely many times
%the order of contact between some roots (at least two) becomes $0$.
%In this way the highest multiplicity of the roots of

We consider  a privileged chart of the blowing-up of the origin in $\R^2$,
more precisely the mapping
$\sigma : (x,y)\mapsto (xy,y)$. Note that by Lemma \ref{orderi} we have
\begin{equation}\label{atilda}
a_i (xy,y) ={\tilde a}_i(x,y)y^i,
\end{equation}
with ${\tilde a}_i(x,y)$ analytic   in a neighbourhood of $0\in \R^2$.
We put
\begin{equation}\label{pptild}\tilde P(x,y,z)= z^d +\sum_{i=1}^{d}\tilde
a_i(x,y)z^{d-i}.
\end{equation}
      We call  $\tilde c_0^i (y)=\frac{c_0^i (y)}{y}$ the { \it proper
transform of}
${c_0^i (y)}$, they
      are the roots of the polynomial
\begin{equation}\label{manyvartilda}
\tilde P(0,y,z)= z^d +\sum_{i=1}^{d}\tilde a_i(0,y)z^{d-i}.
\end{equation}
Clearly
     the  orders of contact between the above  proper transforms drop by 1.
Either all $\tilde c_0^i (y)$ take the same  value at the origin, so by a
shift we may assume that they
vanish at $0$ and we continue our procedure, or
  the highest multiplicity of the roots of the  polynomial
$\tilde P(0,0,z)$ decreases.
Note that the above  procedure has to be finite, as there are at least
  two distinct roots $c_0^i \ne c_0^j$.

%  Let us fix
% one $c_0^i (y)$.
% After finitely many   (say $N_i$) operations the order of contact between
%proper transforms
% $c_0^i (y)$ and $c_0^j (y), \, j\ne i$
% becomes
% $0$ and we are  in the Case 0.
%  More precisely, by Implicit Function Theorem there exists an analytic
%function $F_i:U\to \R$,
% where $U$ is a neighborhood of $0\in \R^2$ such that the proper transform
%of $c_0^i (y)$
% is equal to $F_i(0,y)$.
%
% Denote  by $\sigma^{N_i}= \sigma\circ\cdots \circ \sigma$ composed $N_i$
%times, then $\sigma^{N_i}
% (U\cap
% \{y>0\})$ contains a set
% $H =\{\vert x\vert <y^{N_i+1}, \,0<y<\delta\}$, for
%  $\delta$ is small enough. Hence  $f_i = F_i\circ (\sigma^{N_i})^{-1}$
% is analytic in $H$.

To conclude the first part of the Proposition \ref{keysplitprop}  we have to
explain the remaining case.
\vskip 0,5cm
\item

Case 1.2. Assume that  $c_0^1 (y) = \cdots =c_0^d(y)$.
  By shifting we may assume that
$c_0^1 (y) = \cdots =c_0^d (y)\equiv0$.
% Performing the
%operation described in Case 1.1. we can separate, by Splitting Lemma, the
%proper transform of
%$c_0^1 (y) $ from the proper transforms  of the other  roots $c_0^j
%(y),\,  p<j\le d$.

%shifting by $ z\mapsto z-c^1_0 (y)$ we
%may assume that
%$c_0^i (y)\equiv 0, i=1,\dots , p$.
%Accordingly we may suppose that our hyperbolic polynomial is of the form
%\begin{equation}\label{manyvartilda1}
  %    Q(x,y,z)= z^p +\sum_{i=1}^{p}b_i(x,y)z^{p-i},
%\end{equation}
%where $b_i (x,y)$ are real analytic in a neighbourhood of $0\in \R^2$.
Note that $a_i(0,y)\equiv0$ for all $i$.
After Tschirnhausen transformation we may assume that $a_1 (x,y) \equiv 0$.
%Note that $b_1 (0,y) \equiv 0$,  since all roots of the polynomial
%$z\mapsto P(0,y,z)$ are equal and their sum is $0$.

Now we consider the  coefficient $a_2$; either $a_2 (x,y) \equiv 0$  and
then, by Lemma \ref{hyproots},
polynomial $P$ has the only root
$z(x,y) \equiv 0$ and we are done, or $a_2 (x,y) \not\equiv 0$. In the
second case there  exists an integer $k$
such that
\begin{equation}\label{stop}y\mapsto \frac{\p^k a_2}{\p x^k}(0,y)
\not\equiv 0,
\end{equation} does not vanish identically. We take
the smallest such an integer. By the equation \eqref{znak} $k$ must be even, so
we write $k=2r$.

Applying Lemma \ref{orderi} to
our polynomial
%$(x,z)\mapsto Q(x,y,z)$,
with
$y$ fixed, we obtain that
     \begin{equation}\label{xtilda}
a_i (x,y) ={\tilde a}_i(x,y)x^{ri},
\end{equation}
with ${\tilde a}_i(x,y)$ analytic   in a neighbourhood of $0\in \R^2$.
%
%hence also $\tilde b_1 (x,y) \equiv 0$. Analogously to the identity
%\eqref{newton} we have also
%\begin{equation}\label{a2}
%-2b_2(x,y) = z_1(x,y)^2+ \cdots +z_p(x,y)^2,
%\end{equation}
% where $z_1(x,y), \cdots, z_p(x,y)$ are the roots of
%So if $b_2(x,y)\equiv 0$, then
%where, by lemma \ref{orderi},
Now we consider the polynomial
\begin{equation}\label{manyvartilda2}
\tilde P(x,y,z)= z^p +\sum_{i=2}^{p}\tilde a_i(x,y)z^{p-i},
\end{equation}
     Note that, by \eqref{stop}, we know that $y\mapsto\tilde a_2 (0,y) \not
\equiv 0$.
As a consequence, by Lemma \ref{hyproots}, the polynomial
$(y,z)\mapsto \tilde P(0,y,z) $ has at least 2 distinct roots
     $\tilde c_0^i (y)$ so we may apply the argument of the Case 1.1 and we
are done by the induction on the highest multiplicity.

So we have proved  the following: there exists an integer $N$ such that
\begin{equation*}P(xy^N,y,z)= z^d
+\sum_{i=1}^{d}a_i(xy^N,y)z^{d-i}= \prod_{i=1}^d [z-g_i (x,y)],
\end{equation*}

where $g_i$ are analytic in a neighbourhood of $0\in \R^2$.

So $f_i(x,y) = g_i(xy^{-N},y)$ are the   functions we claimed in
Proposition \ref{keysplitprop}.

We have a more precise control of the functions $f_i$.
Replacing $f_i$ by $f_i (x,y) - f_i(0,y)$ we may assume that
$f_i(0,y)\equiv 0$. (Note that $(x,y) \mapsto f_i(0,y)$ is analytic in
a neighbourhood of $0\in \R^2$.)
Now we can define a strict transform of $f_i$ as
$$f_i^{(1)}(x,y) = y^{-1}f_i(xy,y).
$$
Observe that $f_i^{(1)}$ is a root of the polynomial $\tilde P$ defined by
\eqref{pptild}.
  We have  again $f_i^{(1)}(0,y)\equiv 0$, so we may define 
$f_i^{(2)}$ a strict
transform of $f_i^{(1)}$, and so on.

So  actually we have proved:

\begin{lemma}\label{strans} For each $f_i$ there exists an integer 
$N$ such that
our   strict transform  $f_i^{(N)}$ is analytic in
a neighbourhood of $0\in \R^2$.
\end{lemma}

%More precisely: let us choose   $(0,y_0)$ which is regular point of the
%discrminant
%$D(x,y)$ of polynomial $P(x,y,z)$. We consider now polynomial
%($x,z),y_0,z)\mapsto P(x,y_0,z)$. By Rellich's theorem we may write
%$$
%P(x,y,z) =\prod_{i=1}^d [z-g_i(x,y)]
%$$
%with $x\mapsto g_i(x,y_0)$ analytic in a neighborhood of $0\in \R$.
%In fact $g_i(x,y)$ analytic in a neighborhood of $(0,y_0)\in \R$.
%
%(Argument ?????
%Puiseux or quasi-ordinary singularities)
%
%Clearly some of $g_i(x,y)$ may be multiple roots so they will not split,
%however
%they will appear as multiple roots $f_i$ in the decomposition
%\eqref{keyspliteq}.

\end{itemize}

\vskip 1cm
{\it Proof of} Proposition \ref{keylipprop}.
Now we shall prove that $y\mapsto \frac{\partial f_i}{\partial x}(0,y)$ is
analytic
at $0\in \R$.
Let us expand $f_i$ as  a power series in $x$,
$$f_i(x,y) = \sum_{n=0}^\infty c^i_{n}(y) x^{n}.$$

We have to prove that $\frac{\partial f_i}{\partial x}(0,y) =c^i_{1}(y)$ is
analytic at $0\in \R$.

By the change of variable $z \mapsto
z-c^i_0(y)$, we may assume that $f_i(0,y)\equiv c^i_{0}(y)\equiv 0$.
  Let us compute a proper transform of $f_i$:
$$\frac{f_i(xy,y)}{y} = \sum_{n=1}^\infty c^i_{n}(y)y^{n-1} x^{n}.$$
So the coefficient $c^i_{1}(y)$  remains unchanged ! We know, by Lemma
\ref{strans}  that our   strict transform  $f_i^{(N)}$ is analytic in
a neighbourhood of $0\in \R^2$.

Thus $c^i_{1}(y)$ is a partial
derivative  of an analytic function $f_i^{(N)}$, hence it is analytic itself.

\vskip 1cm

\section{Roots of hyperbolic polynomials are lipschitz}

%We first state the main result of our paper.
We  answer positively a
question
asked by S.{\L}ojasiewicz. First we introduce   some notations.
Consider a polynomial
\begin{equation*}P(x,z)= z^d +\sum_{i=1}^{d}a_i(x)z^{d-i},
\end{equation*} with $a_i:\Omega\to \R$
real analytic functions in an open set $\Omega\subset \R^n$. Assume that for
each $x\in \Omega$ all the roots of the polynomial $z\mapsto P(x,z)$
are real; we
denote them by
$\lambda_1(x)\le\dots \le \lambda_d(x)$. So we have a  mapping $\Lambda:
\Omega \to \R^d$
defined by
\begin{equation}\label{eigenv}\Lambda(x)=(\lambda_1(x),\dots,\lambda_d(x)).
\end{equation}
By a classical result    we know  that
$\Lambda$ is continuous (see eg. \cite{br},
\cite{loj}).
But of course $\Lambda$ is not analytic;  take for instance $z^2-x^2$, then
$\lambda_1(x) =-|x|$ and $\lambda_2(x) =|x|$. If $n=1$, then by
Rellich's Theorem
\ref{rellich} we
can write the components of $\Lambda $ as a MinMax of a family of $d$ analytic
functions.
But this no longer possible if $n\ge 2$, consider $z^2-(x_1^2+x_2^2)$.
However this example suggests
that
$\Lambda$ is more then merely continuous and S.{\L}ojasiewicz asked whether
     $\Lambda$ is
locally lipschitz. Indeed this  is  the case. We will prove the following.

\begin{theorem}\label{liphip} The mapping $\Lambda:\Omega \to \R^d$ is
locally lipschitz.

\end{theorem}
This result is quite delicate as shown by  several examples of
arc-analytic functions
which are not lipschitz (see Section \ref{arcan}). The proof of the
theorem  will be given in
the next section.
We now relate our theorem to some known  facts in the literature.

\subsection{Lidskii's theorem; hyperbolic polynomials versus
symmetric matrices.}
Let $\mathcal {S}_d$ denote   the space  of $d\times d$ symmetric
matrices  with
real coefficients.
Recall  that $\dim \mathcal{S}_d = \frac{d(d+1)}{2}$. We have a
canonical analytic map
$$\theta: \mathcal {S}_d \to \mathcal P_d$$ which
associate
     to  a matrix $A \in \mathcal {S}_d$   its characteristic
polynomial $\theta(A)$. Here $\mathcal P_d$ stands for the space of
monic polynomials of degree $d$.
We identify a vector in $\R^d$ with a monic polynomial of degree $d$ as in
Section \ref{ssplitting}.
Let us  denote by  $\mathcal H_d =\theta(\mathcal {S}_d) $ the space
of hyperbolic
polynomials.
$\mathcal H_d$ is semialgebraic and can be explicitly described by
inequalities involving
subresultants
\cite{br}.
Actually, $\mathcal H_d$  is the closure of a connected component of
the complement of the
discriminant.  Geometry of   $\mathcal H_d$ was studied by Arnold
\cite{ar}, Givental
\cite{g}, Kostov \cite{ko} and others. Its boundary is
concave and piecewise
lipschitz. We have another  canonical map $$\bar\Lambda : \mathcal
{S}_d \to \R^d,$$ which
associate
     to  a matrix $A \in \mathcal {S}_d$  its eigenvalues in the
increasing order,
as in \eqref{eigenv}.
There is a classical result, known as Lidskii's Theorem,  which asserts
the following.
\begin{theorem}\label{lidskiithm}{\bf (Lidskii 1950)}
Given two symmetric matrices $A,B \in \mathcal {S}_d$ then
$$
(\bar\Lambda(A)-\bar\Lambda (B)) \subset\text{conv}\{\tau(\bar\Lambda(A-B)): \,
\tau \in \mathcal {B}_d\},
$$
where $\mathcal {B}_d$ stands for the group of permutations of the 
$d$ coordinates and
$``conv"$ for the convex hull.
\end{theorem}

     Lidskii's Theorem is not
trivial at all, for proofs see \cite{kato},\cite{bhatia}. In
particular it implies the following.
\begin{corollary}\label{lipskicor}
The mapping $\bar\Lambda : \mathcal {S}_d \to \R^d$ is  globally
lipshitz, (with
an explicit constant).
\end{corollary}
Note that our Theorem \ref{liphip} implies that $\bar\Lambda :
\mathcal {S}_d \to
\R^d$ is  locally lipschitz. Indeed we can write
$\bar\Lambda = \Lambda \circ \theta$, in other words we can consider the
analytic (in fact polynomial) family of characteristic polynomials
parametrized by all symmetric matrices. However note that Lidskii's
theorem does not imply our Theorem. Actually there are analytic families
of hyperbolic polynomials which are not associated to an analytic family
of symmetric matrices. More precisely if $P:\Omega \to \mathcal P_d$ is
an analytic mapping then, in general, there is no analytic mapping
$A:\Omega \to \mathcal {S}_d$ such that $P(x)$ is the characteristic polynomial
of $A(x)$ for any $x\in \Omega$. Of course this is true if $\Omega \subset \R$;
by Rellich's theorem, it is enough to take as  $A(x)$ a diagonal matrix
with the roots of $P(x)$ on diagonal.

F. Rellich  observed in his book
\cite{rellichbook}, Chapter I, Section 2,   the following: let
$a_2$ be an analytic function, then the polynomial
$$P(x,z)=z^2 - a_2(x)$$
is hyperbolic if $a_2(x)\ge0$.
    Assume that $P(x,z)$ is a
characteristic polynomial of an analytic
family of matrices
$$
\begin{pmatrix}
&a(x)& \; &b(x)& \cr
&b(x)& \; &-a(x)& \cr
\end{pmatrix}.
$$
It follows that $a_2(x)=a(x)^2+b(x)^2$.
Rellich proved that any positive analytic function in 2 variables
is a sum of 2 squares of analytic functions. But he also showed
that in  general a positive analytic function in $3$ variables
is not  a  sum of  2 squares of analytic functions.
    This  is related to the Hilbert's
17th problem.

\begin{rema}
Our Theorem
gives a locally lipschitz section $\lambda :\mathcal {H}_d \to
\mathcal {S}_d$ of $$\theta: \mathcal
{S}_d \to \mathcal H_d.$$

%There are counterexamples for hyperbolic polynomials with only smooth
%coefficients (\cite{a} page 6) so our theorem cannot be extended to all
%hyperbolic polynomials and so there is no section from all hyperbolic
%to $\mathcal {S}_d$.
\end{rema}

\section{Proof of theorem \ref{liphip}}

We show that Theorem \ref{liphip} follows from Proposition \ref{keylipprop}.

We will show that the components $\lambda_i(x)$ of  $\Lambda:\Omega\to \R^d$
are locally
lipschitz. That is; for any point $x_0\in \Omega\subset \R^n$ there exists
$r=r(x_0)>0$
and $L = L(x_0)<\infty$ such that, if $|x-x_0|<r$ and $|y-x_0|< r$, then
\begin{equation}\label{lipestim}
|\lambda_i(x)-\lambda_i(y)| \le L|x-y|, \, i =1,\dots, d.
\end{equation}
Recall that $\lambda_i$ are  $C^1$  in $\Omega$, except a nowhere
dense, analytic
subset set $A\subset\Omega$.
We are going to prove that each $\frac{\p \lambda_i}{\p x_k}$ is bounded in a
neighbourhood of
$x_0$, more precisely at points where $ \lambda_i$ is $C^1$, that is
outside the analytic set
$A$. Assume that this not the case for $\frac{\p \lambda_i}{\p x_1}$. 
Note that
$\lambda_i$ has semi-analytic graph, and
then by the curve selection lemma   (conform for instance \cite{BM88},
\cite{L}) it follows
that
there exists an analytic arc  $\gamma:(-\varepsilon,\varepsilon)\to \R^n$ such
that
$$
\gamma(0) =x_0, \, |\frac{\p \lambda_i}{\p x_1}(\gamma(s))|\to \infty,\,
\ as \, s\to 0.
$$
Let $e_1 =(1,0,\dots,0)$, and consider the mapping $g(s,t) =\gamma(s) 
+te_1$ and
the associated hyperbolic polynomial $(s,t,z)\mapsto Q(s,t,z) =P(g(s,t),z)$.
According to the
Proposition \ref{keysplitprop},  it splits in a horned neighbourhood
of $s$-axis
into
$\prod_{i=1}^d[z-g_i(s,t)]$ with
$g_i$ analytic in that neighbourhood. So by Proposition \ref{keylipprop}

$$
\frac{\p \lambda_i}{\p x_1}(\gamma(s)) =\frac{\p g_i}{\p t}(s,0)
$$
is bounded for $s\to 0$. This is a contradiction, hence Theorem \ref{liphip}
follows.

\section{Roots of hyperbolic polynomials as multivalued arc-analytic functions}

In this section we prove that the roots of hyperbolic polynomials can
be desingularized
by sequences of blowing-ups with smooth centers. First we recall some
known facts from algebra.

    \subsection{Generalized discriminants}
Consider a generic polynomial
$$
P_c (z)= z^d +c_1 z^{d-1} + \cdots +c_d
$$
for  $z\in \C$ \and  $c=(c_1 ,\dots ,c_d )\in\C^d$. We put
$$
W_s = \{
c\in \C^d:\,
P_c (z) \text{ has at most $s$   distincts roots }\}.
$$
Let $K=\{1,\dots ,d\}$ and put
$$
\mathcal D_s (z_1 ,\dots ,z_d )=
\sum_{  J\subset K;\\ \# J =d-s }
\,
\prod_{ \mu ,\nu \in J ;\mu<\nu \\  }
(z_\mu -z_\nu )^2
\qquad ,\quad s=0,\dots ,d-1.
$$
Since $\mathcal D_s (z_1 ,\dots z_d )$ is a symmetric polynomial, we have
    $\mathcal D_s = D_s \circ \sigma$ with $\sigma=(\sigma_1 ,\dots ,\sigma_d
)$, where
$\sigma_1 ,\dots ,\sigma_d $ are the basic symmetric polynomials
(by the well known theorem on symmetric functions).
So $D_s$ is a polynomial in $c=(c_1 ,\dots ,c_d )$.
    We shall call the sequence  $D_s(c)$,
    $s=0,\dots ,d-1$  the {\it generalized discrminants} of the 
polynomial $P_c$.
By a  similar theory of subresultants (see eg. \cite{br}) we can find an
explicit expression for $D_s(c)$ as a minor of the Sylvester matrix of
$P_c$ and $P_c^\prime$.
Note that
$D_0(c)$ is the discriminant of
$P_c$.
\begin{lemma}\label{gdlem1} For $s=0,\dots ,d-1$ we have
$$
W_s = \{
c\in \C^d :
D_0 (c) = \cdots = D_{d-s-1} (c) =0
\}.
$$
\end{lemma}

Indeed, if  $c\in W_s$ and  $ z =(z_1 ,\dots ,z_d )$ is the complete
sequence of roots of  $P_c (z)$,  then
$
\# \{z_1 ,\dots ,z_d \} \leq s
$, hence $
\mathcal D_0 ( z ) = \cdots \mathcal D_{d-s-1} ( z ) =0$;
which implies
$$ D_0 (c)=\cdots =D_{d-s-1} (c)=0.$$

Conversely, let  $c\in \C^{ d}$ be such that $D_0 (c) =\cdots =D_{d-s-1}
(c)=0$ and let
     $ z
=(z_1 ,\dots ,z_d )$ the complete
sequence of roots of  $P_c (z)$ .
Assume that $c\notin W_s$, hence  $s+1 \leq \# \{ z_1 ,\dots ,z_d \} =t$.
Let
$z_1 ,\dots z_t$  be the distinct roots of  $P_c (z)$. So
$$
\mathcal D_j (z_1 ,\dots ,z_d ) =D_j (c) =0 \qquad \text{for} \qquad
j=0,1,\dots d-s-1.
$$
Since  $d-t\leq d-s-1$,
$$
0=\mathcal D_{d-t} (z_1 ,\dots ,z_d )=\prod_{ \mu <\nu ,\\ \mu ,\nu
\in
\{ 1,\dots ,t\} } (z_\mu -z_\nu )^2    \, ,
$$
which is a contradiction. By  the same argument we obtain the following.

\begin{corollary}\label{gd2} Assume that $P_c$ has exactly $s$ distinct
roots $z_1 ,\dots z_s$. Denote by
$$\tilde P_c(z)= \prod_{i=1}^s (z-z_i)$$
   a square-free polynomial which has
the same roots as $P_c$, and by  $D\tilde P_c$
   the discriminant of
$\tilde P_c$. Then
   $$\nu_1\cdots \nu_s D\tilde P_c= D_{d-s}(c),$$
where each $\nu_i$ is the multiplicity of  $z_i$ as a root of $P_c$.
\end{corollary}
In particular  we can check whether $D\tilde P_c\ne 0$ without computing the
coefficients of  $\tilde P_c$.
\subsection{Splitting according to multiplicities of roots}
Consider polynomials of the form
\begin{equation*}P(x,z)= z^d +\sum_{i=1}^{d}a_i(x)z^{d-i},
\end{equation*} with $a_i$
holomorphic functions  in an open connected subset $U$ of $\C^n$
(or more generally in a connected holomorphic manifold $U$). Recall that
$\mathcal M (U)$,  the ring of meromorphic functions on $U$ is actually a
field. So in the ring of polynomials $\mathcal M (U)[z]$ we have well
defined $gcd$ (greatest common divisor) of any finite family
    of polynomials in  $\mathcal M (U)[z]$. In
particular, if $P,Q \in \mathcal M (U)[z]$ are monic polynomials with
holomorphic coefficients, then $R=gcd(P,Q)$ is again a monic polynomial
with holomorphic coefficients. Indeed if we assume that $R$ is monic (and
$\deg R\ge 1$) then a priori the coefficients of $R$ are only meromorphic,
but the zeros of $R$ are contained in zeros of $P$ which are locally
bounded (as multivalued functions of $x$). So the coefficients of $R$,
being bounded and meromorphic, are actually holomorphic.

We shall say that $P$ is {\it square-free} if its discriminant
$$
D P(x) = D_0(a_1(x), \dots, a_d(x))\not\equiv 0.
$$
Recall that $DP:U\to \C$ is a holomorphic function.
Of course each polynomial in $\mathcal M (U)[z]$ has a unique (up to
a permutation) decomposition into irreducible factors.

%We have the
%following characterization of irreducible polynomials.
%\begin{lemma}\label{irreducible}
%?????
%\end{lemma}

   We shall need the
following splitting.
\begin{proposition}\label{multisplit}Let $U$ be   an open connected subset
   of
$\C^n$ (or more generally  a connected holomorphic manifold).
Let
\begin{equation*}P(x,z)= z^d +\sum_{i=1}^{d}a_i(x)z^{d-i},
\end{equation*} be a polynomial  with $a_i$
holomorphic in $U$.

   Then there are unique  (up to permutation) square-free
monic polynomials $P_1,\dots, P_k$ with coefficients holomorphic in $U$
and pairwise distinct   integers $\nu_1,\dots,\nu_k\ge 1$,
such that
\begin{equation}\label{multispliteq}P= P_1^{\nu_1}\cdots P_k^{\nu_k}.
\end{equation}
Moreover $P_1,\dots, P_k$ are relatively prime; that is if $i\ne j$ then
$gcd(P_i,P_j)=1$.
\end{proposition}
\
\begin{proof} Let $P' = \frac {\p P}{\p z}$. If $P$ is not square-free then
$DP =0$ in
$\mathcal M (U)$,  so $R=gcd(P,P')$ is of degree at least $1$. Hence
$P= RQ$, where $Q$ is a monic polynomial with holomorphic coefficients
in $U$. But  $\deg R<d$ and $\deg Q<d$, so it is easy to conclude applying
   induction  on degree to $R$ and $Q$. Alternatively we can decompose
\begin{equation}\label{irredecomp}P=Q_1^{m_1}\cdots Q_l^{m_l},
\end{equation}
where $Q_j$ are irreducible. Now for a fixed integer
$\nu_i\in \{m_1,\dots,m_l\}$  we write $P_i$ as the product of those
irreducible factors of $P$ which appear in \eqref{irredecomp} with
exponent $\nu_i$.
\end{proof}

We will denote by $\tilde P = P_1\cdots P_k$ the associate square-free
polynomial.
Of course we have also $\tilde P = Q_1\cdots Q_l$.
   Clearly $P^{-1}(0) =\tilde P^{-1}(0)$. It follows from
Corollary \ref{gd2} that we can compute the discriminant $D\tilde P$
without performing the splitting \eqref{multispliteq}. Precisely,

\begin{corollary}\label{multisplitcor} Assume that $P(x,z)$ is as in
Proposition \ref{multisplit}. Let $s=\sum_{i=1}^k \deg P_i$. Then
$$ \nu_1\cdots \nu_s D\tilde P (x) = D_{d-s}(a_1(x), \dots,a_d(x))
\not\equiv 0.
$$
Moreover for each $x\in U$ the polynomial $z\mapsto P(x,z)$ has at most
$s$ distinct roots  and if $D\tilde P(x)\ne 0$, then it has exactly
$s$ distinct roots.
\end{corollary}

    \subsection{Quasi-ordinary singularities}
Let $U$ be   an open subset
   of
$ \C^n$ (or more generally  a holomorphic manifold).
Let
\begin{equation*}P(x,z)= z^d +\sum_{i=1}^{d}a_i(x)z^{d-i},
\end{equation*} be a  polynomial  with $a_i$
holomorphic in $U$. We say that $P$ is {\it quasi-ordinary}
if
   the discriminant $D\tilde P$ of the  square-free reduction $\tilde P$ of
$P$
is  a normal crossing. In other words for each
$a\in U$ there exists a local chart around  $a$ such that
$D\tilde P(x)= u(x) x_1^{\alpha_1} \cdots x_n^{\alpha_n}$, with $u(a) \ne
0$.

The concept of quasi-ordinary singularities goes back (at least) to Jung's
(1908) desingularization of embedded algebraic surfaces. In fact they
appear as ``terminal" singularities which can be resolved by the
normalization. We shall need a crucial property which generalizes the
Newton-Puiseux parametrization. The result below is sometimes called
Abhyankar-Jung theorem.
%\vskip 1cm
\begin{theorem}\label{jungthm}{ \bf (Jung 1908)}
Let $U = \{|x_1|<r_1\}\times \dots \times \{|x_n|<r_n\}$ be
   an open polydisc  in $ \C^n$ and let
\begin{equation*}P(x,z)= z^d +\sum_{i=1}^{d}a_i(x)z^{d-i},
\end{equation*} be a  polynomial  with $a_i$
holomorphic in $U$. Assume that
   the discriminant $D\tilde P$ of the  square-free reduction $\tilde P$ of
$P$ is of the form
$D\tilde P(x)= u(x) x_1^{\alpha_1} \cdots x_n^{\alpha_n}$, where $u(x)$
is a holomorphic non-vanishing function in $U$. Then there exist  integers
$q_1,\dots, q_n \ge 1$ and  holomorphic
functions
$f_1,\dots, f_d$ defined in the polydisc
$U^\prime = \{|z_1|<r_1^{1/q_1}\}\times \dots \times
\{|z_n|<r_n^{1/q_n}\}$ such that
\begin{equation}\label{jungeq}P(x_1^{q_1},\dots, x_1^{q_n},z)=
\prod_{i=1}^d (z-f_i(x_1,\dots,x_n))
\end{equation}
for any $(x_1,\dots,x_n)\in U^\prime$, $z\in \C$.
\end{theorem}
See for instance \cite{zariski}, \cite{lipman} or
\cite{barth}.

\subsection{Splitting of quasi-ordinary hyperbolic polynomials}

We formulate now an important  consequence of Jung's theorem for
hyperbolic polynomials.

\begin{proposition}\label{jungprop}
Let $\Omega =  (-r,r)^n $ be
   an open cube  in $ \R^n$ and let
\begin{equation*}P(x,z)= z^d +\sum_{i=1}^{d}a_i(x)z^{d-i},
\end{equation*} be a hyperbolic polynomial  with $a_i$
analytic  in $\Omega$. Assume that
   the discriminant $D\tilde P$ of the  square-free reduction $\tilde P$ of
$P$ is of the form
$D\tilde P(x)= u(x) x_1^{\alpha_1} \cdots x_n^{\alpha_n}$, where $u(x)$
is analytic and non vanishing in $\Omega$. Then there exist  analytic
functions
$f_1,\dots, f_d:\Omega\to \R$  such that
\begin{equation}\label{jungeq2}P(x,z)=
\prod_{i=1}^d (z-f_i(x))
\end{equation}
for any $x\in \Omega$, $z\in \C$.
\end{proposition}
\begin{proof} Note that it is enough to prove the result  in a
neighbourhood  of any point of $\Omega$ and then use the uniqueness of
analytic extension to obtain functions $f_1,\dots, f_d$ defined in
$\Omega$. So we may assume that the coefficients $a_i$ are actually
holomorphic in the polydisc
$U = \{|x_1|<r\}\times \dots \times \{|x_n|<r\}$ in $\C^n$. We may also
assume that the unit $u(x)$ does not vanish in $U$. We  apply Jung's
Theorem
\ref{jungthm}. Let us take the smallest integers $q_1,\dots, q_n \ge 1$ so that
there are analytic functions $f_1,\dots, f_n$ such  that
formula \eqref{jungeq} holds. We claim that actually
$q_1=\dots= q_n = 1$.

Assume  that one $q_i>1$, for instance that $q_1>1$. We expand
$f_i$ as a power series in $x_1$ with coefficients holomorphic
in $x'=(x_2,\dots,x_n)$. Let us fix an $i$ and  write $f$ instead of
$f_i$. So we have
\begin{equation}\label{jungexp}
f(x_1,x')=
\sum_{\nu=0}^\infty c_\nu(x')x_1^\nu
\end{equation}

Since $q_1>1$ is minimal, then there exists
$\nu_0 \in \NN\setminus q_1\NN$ such
that $c_{\nu_0}\not\equiv 0$.
So there   exists $a'\in
(-r,r)^{n-1}\subset\R^{n-1}$, such that $c_{\nu_0}(a')\ne 0$.
Hence the  Puiseux expansion  of the function
$g(x_1) = f(x_1^{\frac{1}{q_1}}, a'), \, x_1>0$ has at least one  monomial with
non-integer exponent. So $g(x_1)$ cannot be extended to an analytic
function in a neighbourhood  of $0\in \R$. But on the other hand, by
Rellich's Theorem \ref{rellich},  the roots of  the polynomial $P(x_1,a',z)$
are analytic functions on $(-r,r)$. Clearly  $g$ must be a
restriction of one of these functions. This is a contradiction.

\end{proof}

As a first  consequence of Proposition \ref{jungprop} observe that
the roots of an analytic family of hyperbolic polynomials can  be chosen
analytically outside a subset of codimension at least two.

\begin{theorem}\label{codim2}Consider a polynomial
\begin{equation*}P(x,z)= z^d +\sum_{i=1}^{d}a_i(x)z^{d-i},
\end{equation*}
where $a_i:\Omega\to \R$ are
real analytic functions in an open set $\Omega\subset \R^n$.
Assume that for
each $x\in \Omega$ all roots of the polynomial $z\mapsto P(x,z)$ are real.

Then, there
exists $\Sigma \subset \Omega$ a  semianalytic closed set of codimension
at least
$2$ such that if
$a\in \Omega \setminus \Sigma$ then there is a neighbourhood $U$ of $a$ and
analytic functions $f_i: U\to \R, \, i =1,\dots, d$ such that
\begin{equation*}P(x,z)=
\prod_{i=1}^d (z-f_i (x)), \end{equation*}
for any $ x \in U$, $z\in \R$.
\end{theorem}
\begin{proof}
Without loss of generality we may assume that $\Omega$ is connected. So
   the discriminant $D\tilde P$ of the  square-free reduction $\tilde P$ of
$P$ is a well defined non vanishing analytic function on $\Omega$.
   We are going to prove that $D\tilde P$ is a normal crossing
outside a closed semianalytic set of codimension at least $2$.

Let $Z$ be the set of zeros of $D\tilde P$. Clearly  $Z$ is an analytic
subset of $\Omega$. Let
$Reg_{n-1} Z$  be  the set  of points  $x\in \Omega$ such that for
some neighbourhood $U$ of $x$ the set  $Z\cap U$ is an analytic submanifold
of dimension $n-1$.  Of course $Reg_{n-1} Z$ is
open in $Z$ and by {\L}ojasiewicz's theorem \cite{L},
$$
\Sigma' = Z\setminus Reg_{n-1} Z,
$$
is a semianalytic set, $\dim \Sigma'< \dim Z \le n-1$. Hence
$\Sigma'$ is closed in  $\Omega$ of codimension at least $2$.
Let $\Delta $ be a connected component of  $Reg_{n-1} Z$.
Let $\alpha $ be the smallest integer such that
$$
h =\frac{\partial^\alpha D\tilde P}{x_1^{r_1}\dots x_n^{r_n}}
$$
does not vanish identically  on $\Delta$ for some multi index
$(r_1, \dots,r_n)$, $\alpha =r_1+ \cdots +r_n$.
Hence $\Sigma''(\Delta) = h^{-1}(0)\cap \Delta$ is a semianalytic set  of
dimension less than $(n-1)$.
Note that if $a\in \Delta \setminus\Sigma''(\Delta)$ then in some chart
around $a$ we can write
$$ D\tilde P(x_1, \dots, x_n) = u(x) x_1^\alpha,
$$
with some unit $u(x)$.
   Let
$$
\Sigma''=\bigcup \Sigma''(\Delta),
$$
where the  union is taken over all  connected components of $Reg_{n-1} Z$.
Finally we put
$$
\Sigma = \Sigma'\cup \Sigma''.
$$
Clearly  $\Sigma$ is  semianalytic and closed in  $\Omega$ of codimension at
least
$2$.  Let $a\in \Omega \setminus \Sigma$,
%, then either $D\tilde P(a)\ne 0$
%and, by implicit function theorem, $P(x,z)$ splits in a neighbourhood  of
%$a$ into a product of linear factors or $D\tilde P(x)$ is of the
%form $u(x)x_1^\alpha$ and
then $D\tilde P$ is a normal crossing in a neighborhood of $a$ so we conclude
using Proposition
\ref{jungprop}.
\end{proof}
\begin{rema}In particular if $\dim M =2$ then $\Sigma$ has only isolated
points. In other words  any
$2$-parameter analytic family of hyperbolic polynomials  splits
   locally, outside a
discrete set, into linear factors. More generally, observe
that if
$\Omega$ is connected  and
$\Sigma$ is  semianalytic closed in  $\Omega$ of codimension at least
$2$, then $\Omega\setminus \Sigma$ is also connected. However we
cannot claim that we can  split $P(x,z)$  into a product of linear factors
in $\Omega\setminus \Sigma$. Here we may have a nontrivial monodromy.
\end{rema}
\begin{exam} The discriminant of the hyperbolic polynomial
$$P(x,y,z)=z^3-3(x^2+y^2)z-2x^3,$$
vanishes on the $x$-axis. Here $\Sigma $ is just the origin.  Note 
that $P(x,0,z)=
(z+x)^2(z-2x)$. So for
$x>0$ the double root is smaller than the simple root, while  for 
$x<0$ their order is
inversed. Moving around  the circle $\{x^2+y^2=1\}$ in $\R^2$ gives a 
nontrivial
monodromy.
\end{exam}

     \subsection{Multivalued arc-analytic functions} For the purpose of studying
hyperbolic
polynomials we use the following notion. Let $M^m, N^n$  be two real
analytic
manifolds.
Let $F$ be a subanalytic subset of $M\times N$. For $x\in M$ we  denote
$$F(x) = \{y\in N:\,  (x,y)\in F\}$$
and  we call $F(x)$ the {\it set of values of $F$ at $x$}. If $F(x)$ is non
empty for every $x\in M$
we say that $F$ is {\it a multivalued mapping on $M$ with values in N}. If
$M$ is connected we say that
$F$ is {\it
$k$-valued} if $F(x)$ has at most $k$ points for any $x\in M$ and exactly $k$
points for some
$x_0\in M$. Single valued $F$ is a function in the usual sense.

We will say that $F$ is {\it continuous} if $F$ is closed in  $M\times N$.
We call
$F$  {\it proper} if the projection on $M$ restricted to $F$ is a proper map.
We say that $F\subset M\times N$ is a {\it  $k$-valued arc-analytic
mapping} if
for any  analytic arc $\gamma:(-\varepsilon, \varepsilon) \to U \subset M$
     there are $k$
analytic functions $f_i:(-\varepsilon, \varepsilon) \to N,\, i=1,\dots,k$
such that
$$F(\gamma(t)) = \{f_1(t),\dots,f_k(t)\}.$$
Note that, in general, the set $\{f_i(t)\}$ is {\bf not ordered}. If $F$
is single  valued then it is an
arc-analytic mapping in the usual sense.

\vskip 1cm
\begin{theorem}\label{multires} Every proper  $k$-valued arc-analytic and
subanalytic mapping
is  locally blow-analytic via a locally finite
collection of analytic modifications $ \sigma_{\alpha}: W_{\alpha} \to
M$, that is for any $\sigma_\alpha$ we have
$$\tilde\sigma_\alpha^{-1}(F) = \bigcup \tilde F_i$$
and each $\tilde F_i$ is a graph of an analytic function in
$\sigma_{\alpha}^{-1}(M)$.
Here $\tilde\sigma_\alpha:W_\alpha \times N \to M\times N$,
$\tilde\sigma_\alpha(w,y) =(\sigma_\alpha(w),y)$. If $F$ is semialgebraic
then, instead of the family $ \sigma_{\alpha}$ we can take one
$\sigma$ which is a finite composition of global blowing-ups with smooth
centers.

\end{theorem}
The proof of this result will be published separately.

In the case of hyperbolic polynomials with analytic coefficients  the Theorem
\ref{multires}
can be restated as follows.

\begin{theorem}\label{hipres}Consider a polynomial
\begin{equation*}P(x,z)= z^d +\sum_{i=1}^{d}a_i(x)z^{d-i},
\end{equation*}
where $a_i:\Omega\to \R$ are
real analytic functions in an open set $\Omega\subset \R^n$.
Assume that for
each $x\in \Omega$ all the roots of the polynomial $z\mapsto P(x,z)$ are real.
Then, there
exists
    $ \sigma: W \to
\Omega$  a locally finite composition of blowing-ups with smooth (global)
centers, such that for any $w_0\in W$ there are a neighbourhood $U$ and
analytic functions $F_i: U\to \R,\, i=1,\dots,d$ such that
\begin{equation*}P_\sigma(w,z)= z^d +\sum_{i=1}^{d}a_i(\sigma (w))
z^{d-i}=
\prod_{i=1}^d (z-F_i (w)), \end{equation*}
for any $ w \in U$, $z\in \R$.
\end{theorem}
\begin{rema}\label{antihipc}\rm Note that the Theorem \ref{hipres}
applies also to real analytic families of monic polynomials
such that for
each $x\in \Omega$ all the roots of the polynomial $z\mapsto P(x,z)$
are  purely
imaginary.  Indeed, if
$P(x,z)$ is such a polynomial then $P(x,iz)$ is hyperbolic with real analytic
coefficients.
\end{rema}

\begin{proof}
Without loss of  generality we may assume that $\Omega$ is connected.
   So
   the discriminant $D\tilde P$ of the  square-free reduction $\tilde P$ of
$P$ is a well defined non vanishing analytic function on $\Omega$.
   By Hironaka's Desingularization Theorem \cite{H},
there exists
    $ \sigma: W \to
\Omega$  a locally finite composition of blowing-ups with smooth (global)
centers, such that $D\tilde P\circ \sigma$ is a normal crossing.
But $D\tilde P\circ \sigma$ is the discriminant of the square-free
reduction of
$P_\sigma(w,z)= z^d +\sum_{i=1}^{d}a_i(\sigma (w))z^{d-i}$. So the
theorem follows immediately from Proposition \ref{jungprop}.
\end{proof}
\begin{exam} \rm Let
$P(z,x_1,x_2) =z^2 -(x_1^2+x_2^2)$, so $ DP= 4 (x_1^2+x_2^2)$ is the
discriminant of  $P$. Clearly the blowing-up of the origin makes it
normal crossing. Namely, we write
   $x_1=w_1$,$x_2=w_1w_2$ for the blowing-up, so
$$
P(w_1,w_2,z)= (z-F_1(w)) (z-F_2(w)),
$$
where $F_1 =w_1(1+w_2^2)^{1/2}$, $F_2 = -w_1(1+w_2^2)^{1/2}$  are real
analytic functions (defined in one chart). Note that these functions are
not holomorphic  if we consider $w_1,w_2$ as complex numbers.
This simple example shows the purely real character of Theorem \ref{hipres}.
\end{exam}
   In the sequel we shall need that the
analytic functions $F_i: U\to \R$ in Theorem \ref{hipres} are also
normal crossings. This can be  achieved by making the coefficient $a_d$
normal crossing, indeed we have $a_d\circ \sigma = F_1\cdots F_d$. Of
course each factor of  a normal crossing is again a normal crossing.
Assume that we arranged $F_i$'s in  such a way that  $F_1, \dots,F_s$   are
all the distinct roots of $P_\sigma(w,z)$.
%Assume  the $m_i$ is the multiplicity of $F_i$, for  $i\le s$.
Recall that,
% (cf. Corlollary \ref{gd2})
$$D\tilde P\circ \sigma (w) =  \prod_{i<j\le s}
(F_i(w)-F_j(w))^2.$$
But $D\tilde P\circ \sigma (w)$ is a normal crossing  so it follows
that for any $  i,j\le s, \,i\ne j$ the function $(F_i(w)-F_j(w))$
is a normal crossing as well. Recall now a very important observation in
Bierstone and Milman \cite{BM88}:
\begin{lemma}\label{wellorder}
Let $U$ be an open connected subset of $\R^n$.
Let  $F_i\not\equiv 0,\,
i=1,\dots,d$  be  analytic functions in $U$.
Assume that all $F_i$ and all their differences $F_i
-F_j$ are normal crossings (or identically $0$). Then for each $w\in U$
the exists a neighbourhood $U_w$ such that
for any $i,j\le d$ at least one of the functions $\frac{F_i}{F_j}$ or
$\frac{F_j}{F_i}$ extends to an analytic function in $U_w$.
In particular there exists $i_w \le d$ such that
$\frac{F_j}{F_{i_w}}$ extends to an analytic function in $U_w$, for any
$j$. We will say for short that $F_1,\dots, F_d$ are
{\bf well ordered on } $U_w$.
\end{lemma}
As consequence we obtain:
\begin{rema}\label{differences} In Theorem \ref{hipres} we can choose
$U$ and the analytic functions $F_i:U\to \R$ in such a way that they are  well
ordered on
$U$.
\end{rema}

\section{Diagonalization of analytic families of symmetric matrices}

The goal of this section is to generalize a result of Rellich
\cite{rellichbook} which states that a $1$-parameter analytic family
of symmetric matrices admits a uniform diagonalization.

We will denote by $\mathcal{S}_d$ the space of symmetric $d\times d $
matrices with real entries.
We consider first an analytic family of  symmetric matrices  $A:\Omega \to
\mathcal {S}_d$, where
$\Omega$ is an open connected subset of $\R^n$. Assume that the eigenvalues
of $A(x)$ can be chosen analytically in  $\Omega$. Precisely we assume
that there are analytic functions $F_i: \Omega \to \R$, $i=1,\dots,d$ such
that $\{F_1(x), \dots,F_d(x)\}$ is the set of the eigenvalues of $A(x), \,x\in
\Omega$. For a generic point  $x\in \Omega$ (i.e. outside a nowhere
dense analytic subset) each $F_i(x)$   is of the same constant
multiplicity $\nu_i$.
If an eigenvalue $F_i(x)$ is actually of a constant multiplicity on
$\Omega$ then, for any $x\in \Omega$
$$ V_i(x) =Ker (A(x)-F_i(x)\I_d)
$$
is an analytic  family of $\nu_i$-dimensional eigenspaces of
$A(x)$. In particular we can choose locally, in an analytic way,
an orthonormal basis of $V_i(x)$. However in general at some points the
dimension of $ V_i(x)$ may be strictly greater than
$\nu_i$. If  $\Omega=I$ is an interval in  $\R$  and $x_0\in I$ is such
that $\dim V_i(x_0)>\nu_i$, then
$\lim_{x\to x_0} V_i(x)$ exists, in the corresponding Grassmanian.
Moreover the mapping $x\mapsto V_i(x)$ obtained by the continuous extension is
actually analytic (as we will show later on).
However this is no longer true if we consider an analytic family depending
on $n\ge 2$ parameters. Indeed we have,
\begin{exam}\label{uncontrolledeigenspaces}\rm
   Consider a family of symmetric matrices of the form
$$
A(x_1,x_2)=\begin{pmatrix}
&x_1^2& \; &x_1x_2& \cr
&x_1x_2& \; &x_2^2& \cr
\end{pmatrix}, \,( x_1,x_2)\in \R^2.
$$
Note that $\phi = 0$, $\psi=x_1^2+x_2^2$  are eigenvalues of
$A(x_1,x_2)$ and $\Phi=(1,\frac{x_2}{x_1})$,$\Psi=(1,-\frac{x_1}{x_2})$
are the corresponding eigenvectors. Clearly there is no   limit of
$\Phi$ and  $\Psi$ as $( x_1,x_2)\to (0,0)$. So this family cannot be
simultaneously diagonalized in an analytic (even continuous) way. However
if we blow up the origin in $\R^2$, that is we put
   $x_1=w_1$,$x_2=w_1w_2$, then the corresponding family
$$
A(w_1,w_2)= w_1^2\begin{pmatrix}
&1& \; &w_2& \cr
&w_2& \; &w_2^2& \cr
\end{pmatrix}, \,( w_1,w_2)\in \R^2,
$$
admits a simultaneous analytic diagonalization.
\end{exam}

The next  theorem explains  that this happens
for a general analytic family of symmetric matrices.

To fix the terminology we recall that if $\Phi:E\to E$ is a linear mapping and
$\lambda$ is an eigenvalue of  $\Phi$, then
$E_\lambda =\{x\in E; \, \Phi(x)=\lambda x\}$ is called {\bf the eigenspace}
of $\Phi$
(associated to $\lambda$). Any nontrivial linear subspace of $E_\lambda$
is called {\bf an eigenspace}  of $\Phi$
(associated to $\lambda$).

\begin{theorem}\label{diagres}
Consider an
analytic family $A:\Omega \to \mathcal {S}_d$ of symmetric matrices,
where  $ \Omega$ is an open
connected subset of $ \R^m$ and
$\mathcal {S}_d$ stands for the space of
symmetric $d\times d$ matrices with real entries.
   Then, there
exists
   $ \sigma: W \to
\Omega$  a locally finite composition of blowing-ups with smooth (global)
centers, such that for any $w_0\in W$ there is a neighbourhood $U$  such
that  the corresponding family
$ A \circ \sigma |_{U}:U \to
\mathcal {S}_d$ admits a simultaneous
analytic diagonalization.

   More precisely, let
$P_\sigma (w,z)$ be the characteristic polynomial of
   $ A \circ \sigma (w)$; recall that for a generic $w\in W$ the polynomial
$P_\sigma (w,z)$ has  $s$ distinct real roots with the constant number of
roots of fixed multiplicity. Then, for each $w\in W$, there exits an
orthogonal decomposition

\begin{equation}\label{espropres}
\R^d = V_1(w) \oplus \dots \oplus V_s(w), \,
\end{equation}
such that:
\begin{enumerate}
\item each $V_i(w)$ is an eigenspace of $ A \circ \sigma (w)$,
$\dim V_i(w) =m_i\ge 1$;
\item  if the eigenvalue $\lambda_i(w)$ associated to $V_i(w)$ is a root
of $P_\sigma (w,z)$ of multiplicity $m_i$ then,
$$ Ker(A\circ\sigma (w)-\lambda_i(w)\I_d)= V_i(w);$$
\item  if the eigenvalue $\lambda_i(w)$ associated to $V_i(w)$ is a root
of $P_\sigma (w,z)$ of multiplicity $>m_i$ then,
$$Ker(A\circ\sigma (w)-\lambda_i(w)\I_d)= V_i(w)\oplus V_{i_1}(w)\oplus
\cdots\oplus V_{i_k}(w)
$$
for some $i_1,\dots,i_k \in \{1,\dots, s\}\setminus \{i\}$;

%\item $V_i= \bigcup_{w \in W}  V_i(w)\times\{w\}$ is a locally trivial
%analytic vector bundle over $W$.
\item for  any $w_0\in W$ there is a neighbourhood $U$ and
analytic functions $e_i:
U \to (\R^d)^{m_i}$, $i=1,\dots,s$
   such that $ e_i(w)$ is an orthonormal basis of $V_i(w)$.

\end{enumerate}
\end{theorem}
\begin{rema}\label{diagresremark}
Theorem \ref{diagres} holds also for real analytic families of
Hermitian matrices.
\end{rema}

\begin{rema}\label{diagresremark1} \rm We can describe a global structure of
bundles given by \eqref{espropres}. Recall that $W$ is connected so
the polynomial $P_\sigma (w,z)$ admits a unique decomposition into
irreducible factors
\begin{equation}\label{irredecomp1}P_\sigma =Q_1^{m_1}\cdots Q_l^{m_l}.
\end{equation}

Let us fix one $Q_j(w,z)$ and write $m=m_j$. By  Theorem \ref{hipres} we
may assume that in an open set $U\subset W$ we can choose roots of
$Q_j(w,z)$ as  analytic functions $\lambda_1,\dots \lambda_{d_j}:U\to
\R$. Recall that for a  generic $w\in U$ all $\lambda_i(w)$ are simple
roots of $Q_j(w,z)$.
Let us denote by $\GG^m_d$ the Grassmanian of $m$-dimensional subspaces
of $\R^d$. Then according  to Theorem \ref{diagres} we have
   orthogonal subspaces
$V_1(w), \dots, V_{d_j}(w) \in \GG^m_d$, $d_j=\deg Q_j$, which are proper
subspaces of $A\circ\sigma (w)$  associated to $\lambda_1 (w),\dots
\lambda_{d_j} (w)$. Let us collect them together and write
   $$
\Xi_j = \bigcup_{w\in W} \{V_1(w), \dots, V_{d_j}(w)\}\times {w}
\subset \GG^m_d\times W.
$$
Note that  $\Xi_j$ is an analytic submanifold of $\GG^m_d\times W$,
moreover the natural  projection $\pi:\Xi_j \to W$ is an analytic
$d_j$-sheeted covering.
Indeed $\Xi_j$ may be seen as a multivalued analytic function;  each
$U\ni w\mapsto V_i(w)$ is analytic, and the values are distinct since
the subspaces are orthogonal.
Finally we observe that,
\end{rema}
\begin{proposition}\label{cnxe}
$\Xi_j$  is connected.
\end{proposition}
\begin{proof}   Indeed, if $\Xi$ is a connected component of $\Xi_j$, then
   $\pi:\Xi \to W$ is again an analytic $p$-sheeted covering. For each
$w_0$ there is a neighbourhood $U$ and $p$ distinct analytic sections
of $\pi:\Xi \to W$. To each  section (of eigenspaces) we can associate
   the corresponding eigenvalues $\lambda_1,\dots \lambda_{p}:U\to
\R$, which are analytic functions. We put
$$
Q(w,z)= \prod_{i=1}^p(z-\lambda_i(w)), \, w\in U.
$$
By connectedness  we can extend analytically $Q$ on $W$. Clearly $Q$
divides $Q_j$ but
$Q_j$ is irreducible, so $Q=Q_j$, hence  $\Xi=\Xi_j$.

\end{proof}

We begin now the proof of Theorem \ref{diagres}.
Since the characteristic polynomial of a symmetric matrix is
hyperbolic, by  Theorem
\ref{hipres}, there exists
   $ \sigma: W \to
\Omega$,  a locally finite composition of blowing-ups with smooth (global)
centers, such that the eigenvalues of  the corresponding family
$ A \circ \sigma$ are  locally  normal crossings.
We are going to construct $ \sigma': W' \to W$,  a locally finite
composition of
blowing-ups with smooth (global) centers, such that  the
corresponding family
$ A \circ (\sigma \circ \sigma')$  admits a simultaneous
analytic diagonalization.

Recall that we have the splitting of $P_\sigma =Q_1^{m_1}\cdots Q_l^{m_l}$
   into irreducible factors, where $P_\sigma$ is the characteristic
polynomial of the
   family $ A \circ \sigma$. Let us fix one $Q_j(w,z)$ and write $m=m_j$.

We shall
first explain the simpler case where $m=1$. Hence, for a
generic $w\in U$, all eigenspaces
associated to the roots of  $Q_j(w,z)$ are of dimension $1$. Let
$\lambda:U\to
\R$ be an analytic choice of roots of
$Q_j(w,z)$, where $U$ is an open neighbourhood of some fixed point
$w_0\in W$. The
eigenspace of
$A\circ
\sigma(w)$ associated to
$\lambda (w)$ is a set of solutions of a $d\times d$ system
\begin{equation}\label{eigenspace1}
(A\circ\sigma (w)-\lambda(w)\I_d)X= 0
\end{equation}
Recall that for a  generic $w\in U$ this system is of  rank $d-1$. So we can
delete  one equation from \eqref{eigenspace1} and we obtain an
equivalent system
\begin{equation}\label{eigenspace2}
B(w)X= 0,
\end{equation}
where $B(w)$ is a matrix with $d-1$ rows and $d$ columns. Let $M_k(w)$ denote
the determinant of the $(d-1)\times (d-1)$  matrix obtained from
$B(w)$ by deleting
the $k$-th column. By Cramer's rule we obtain that
$$
\bar v(w)= (-M_1(w), \dots ,(-1)^k M_k(w),\dots ,(-1)^d M_d(w))
$$
is  a solution of \eqref{eigenspace1}.
%the eigenspace of $A\circ\sigma(w)$ associated to $\lambda (w)$.
But of course we have to
check that $\bar v(w)\ne 0$, which  is true for a generic $w\in U$,
but in general
not for all $w\in U$. In particular we might possibly have $\bar
v(w_0)= 0$. So we want to divide
all the coefficients of   $\bar v(w)$ by one of them and get again
analytic coefficients.

   We may assume (we
explain it below)  that all minors  $M_k$, $k=1,\dots,d$ are normal
crossings and
moreover that they are well ordered at $w_0$ (cf. Lemma \ref{wellorder}).
Permuting, if necessary, the coordinates in $\R^d$ we may assume that
$M_1(w)$ is the smallest among all  $M_k$, $k=1,\dots,d$. In other words
$$
m_k(w)=(-1)^{k-1}\frac{M_k(w)}{M_1(w)}, \quad k=2,\dots,d
$$
extend  to  analytic functions in a neighbourhood of $w_0$. Thus
$$
    v(w)= (1,m_2(w),\dots ,  m_d(w))
$$
is actually an eigenvector of  $A\circ
\sigma(w)$ associated to
$\lambda (w)$. Clearly $v(w)$ is analytic in a neighbourhood of $w_0$.
Finally  we normalize $v(w)$ in order to get an orthonormal basis of
$V(w) =\R v(w)$, the subspace generated by $v(w)$. Note that, for a
generic $w\in U$, we have
$V(w)=Ker (A\circ\sigma (w)-\lambda(w)\I_d)$.

We consider now the general case where the factor $Q_j$ appears with exponent
$m=m_j\ge1$. So now, for a
generic $w\in U$, $\lambda(w)$ is a root of $Q_j(w,z)$ of multiplicity $m$.
   The eigenspace of
$A\circ
\sigma(w)$ associated to
$\lambda (w)$ is a set of solutions of the $d\times d$ system
\begin{equation}\label{eigenspace3}
(A\circ\sigma (w)-\lambda(w)\I_d)X= 0
\end{equation}
and  for a  generic $w\in U$ this system is of  rank $d-m$. We shall construct
   linearly independent
$v_1(w),\dots,v_m(w)\in \R^d$ which are  analytic in a
neighbourhood of
$w_0$ and
   such that, for a generic $w$ we have
$$
V(w) = Ker (A\circ\sigma (w)-\lambda(w)\I_d),
$$
where $V(w)= span(v_1(w),\dots,v_m(w))$ stands for the subspace  generated
by $v_1(w),\dots,v_m(w)$.

   So we can
delete  $m$ equations from \eqref{eigenspace3} and we obtain an equivalent
(generically) system
\begin{equation}\label{eigenspace4}
B(w)X= 0,
\end{equation}
where $B(w)$ is a matrix with $d-m $ rows and $d$ columns. As in the case $m=1$
we can consider all $(d-m)\times (d-m)$ minors of $B(w)$ and we may assume that
they are well ordered at $w_0$. Let $M(w)$ be the smallest (at $w_0$)
among all these
minors. Permuting, if necessary, the coordinates in $\R^d$ we may suppose that
$M(w)$ is the determinant of the matrix $C(w)$ formed by the first
$d-m$ columns. We
construct a vector
\begin{equation}\label{eigenspace5}
v_1(w) = (a_1(w),1,0,\dots,0),
\end{equation}
where $a_1(w)\in \R^{d-m}$ is the solution of a system
\begin{equation}\label{eigenspace6}
C(w)X'=b_{d-m+1}(w),
\end{equation}
here $b_{d-m+1}(w)$ denotes the $({d-m+1})$-column of $B(w)$. Observe that the
coordinates of $a_1(w)$ are quotients of some minors of $B(w)$ by $M(w)$,
   so they extend to analytic functions in a neighbourhood of $w_0$. 
We construct
$v_2(w),\dots,v_m(w)$ analogously by shifting $1$ to the right in
\eqref{eigenspace5} and considering next columns of $B(w)$ in
\eqref{eigenspace6}.
Finally to obtain an orthonormal basis of $V(w)$ we apply the Gram-Schmidt
orthonormalization to the family $v_1(w),\dots,v_m(w)$.

We are left with proving that after a suitable composition of blowing-ups with
smooth global centers, the  minors of $B(w)$ are normal crossings
which are well
ordered  at any point of $W$.
We begin with a lemma which is actually a description of the normaliztion of a
zero set of an irreducible hyperbolic polynomial with discriminant
which is a normal
crossing.
\begin{lemma}\label{normalization1} Consider a hyperbolic polynomial
$Q(x,z)= z^d +\sum_{i=1}^{d}a_i(x)z^{d-i},$
where $a_i:\Omega\to \R$ are
real analytic functions in a connected analytic manifold $\Omega$.
Assume that $Q$ is irreducible   and moreover that the discriminant 
$DQ:\Omega\to \R$
is a normal crossing. Then there exist a connected analytic manifold 
$\Xi$ and an
analytic
$d$-sheeted covering $p:\Xi\to \Omega$, an analytic function
$z:\Xi\to \R$ such that
$$
Q(p(\xi),z (\xi))=0, \quad \xi \in \Xi.
$$
\end{lemma}
\begin{proof} We define
  $\Xi$ as a space of germs $f_x $, at points  $x\in \Omega$, of analytic
functions
$f:U\to \R$ such that $Q(x,f(x))=0,\, x\in U$, where $U$ are   open subsets of
$\Omega$. We have also  a canonical map $F:U\ni x\mapsto f_x \ni \Xi$,
where $f_x$ stands for the germ of $f$ at the point $x$. These maps
$F$ define an
analytic atlas
on $\Xi$, thus we obtain
a structure of an analytic manifold on $\Xi$. Note that we did not
specify the
topology on $\Xi$, but actually this is not necessary, see for
instance \cite{loj}. We
only have to check that the topology we obtain is Hausdorff, but this
is the case
since we consider only analytic functions.

   Now the mapping
$p:\Xi\to\Omega$ is defined, in the above chart, as the inverse of
$F$, so clearly it is a local diffeomorphism. We put $z(\xi) =f(x)$, for $\xi
=f_x$. It follows from Proposition \ref{jungprop} that
$p:\Xi\to\Omega$ is  indeed
a $d$-sheeted covering. To prove that $\Xi$ is connected we may use
the same argument
as in the proof of Proposition \ref{cnxe}.
\end{proof}

We come back to the proof of  the fact that the minors of the
matrices considered in
\eqref{eigenspace2} and  \eqref{eigenspace4} can be made (by a composition of
suitable blowing-ups) well  ordered normal crossings. Note that these
minors can be
seen as
$d$-valued analytic
   functions on $W$. More precisely they extend to  analytic functions
$M_k:\Xi_j\to \R$  on
the space $\Xi_j$ associated to the polynomial $Q_j$, by Lemma
\ref{normalization1}.
We denote by $p_j:\Xi\to W$ the corresponding covering.
We will  consider only those $M_k,\, k= 1,\dots, K$ which are non
identically zero.
Recall that $\Xi_j$ is connected hence these minors are non
identically zero on any
open subset of $\Xi_j$. So now we can associate to $Q_j$ and $A\circ\sigma$ two
analytic  non identically zero functions $\Phi_j,\Psi_j:W\to\R$
defined as follows
\begin{equation}\label{minor1}
\Phi_j(w) = \prod_{k=1}^K\prod_{\xi\in p_j^{-1}(w)}M_k(\xi)
\end{equation}
and the function which is the product of differences of all factors
in \eqref{minor1},
that is
\begin{equation}\label{minor2}
\Psi_j(w) = \prod (M_{k}(\xi)-M_{k'}(\xi')),
\end{equation}
where the product is taken over all $k,k' \in \{1,\dots, K\}$, $k\ne k'$ and
$\xi,\xi'\in p_j^{-1}(w)$, $\xi\ne\xi'$.
Finally we can take $\Phi$ and $\Psi$ which are respectively the
products of all
$\Phi_j$ and $\Psi_j$ associated  to the prime factors of $P_\sigma$.

By Hironaka's desingularization theorem there exists  $\sigma':W'\to
W$ which is a
locally finite composition of blowing-ups with smooth global centers
such that both
$\Phi\circ\sigma'$ and $\Psi\circ\sigma'$ are normal crossings. Hence in
particular each factor  in \eqref{minor1} and \eqref{minor2} becomes a normal
crossing.

Thus we achieved a proof of Theorem \ref{diagres}.

\begin{rema}\label{diagresremark2} \rm Note that in the proof of Theorem
\ref{diagres} we used not only the fact that the characteristic
polynomials of symmetric
matrices are hyperbolic but also the  fact that the eigenspaces
associated to different
eigenvalues are orthogonal. Indeed we need to know that subspaces
$V(w)$ and $V'(w)$
which are associated to generically different eigenvalues $\lambda (w)$
and $\lambda'(w)$ are orthogonal also for those $w$ for which
   $\lambda (w)=\lambda'(w)$. This is the case  by   continuity. However if we
consider an analytic family of matrices which are diagonalizable over
reals (but not
symmetric), then it may happen that the subspaces $V(w)$ and $V'(w)$
which have trivial
intersection for generic $w$ may have nontrivial intersection for
some $w_0$ such that
$\lambda (w_0)=\lambda'(w_0)$. So Theorem \ref{diagres} does not
apply to such a
family.
\end{rema}

\begin{exam}\label{nonsymdiag}
The  following one parameter family of diagonalizable matrices
$$
\begin{pmatrix}
&1-x^2& \; &x& \cr
&0& \; &1+x^2& \cr
\end{pmatrix}, \, x\in \R.
$$
have eigenvectors $v_1(x) = (1,0)$, $v_2(x) = (1,2x)$, which form a 
basis of $\R^2$
except for $x=0$. So we cannot choose a basis of eigenvectors in a 
continuous way.

\end{exam}
\vskip 1cm

\section {Reduction of analytic families of antisymmetric
matrices}

The method of diagonalization of analytic families of symmetric
matrices described in
the previous section  applies as well to analytic families of
antisymmetric matrices.
Indeed, the characteristic polynomial of an antisymmetric matrix has
purely imaginary
roots and it is easy to  see (cf. Remark \ref{antihipc}) that Theorem
\ref{hipres}
applies also to real analytic families of monic polynomials with
purely imaginary
roots.

First we recall  briefly some basic facts about antisymmetric matrices.
\begin{lemma}\label{antisym}
Let $A$ be an antisymmetric $d\times d$ matrix with real coefficients, then:
\begin{enumerate}
\item The eigenvalues of $A$ are purely imaginary, moreover if $\mu$
is an eigenvalue of
$A$ then $-\mu$ is also an eigenvalue of $A$.
\item The eigenspaces associated to distinct eigenvalues are orthogonal.
\item There exists an orthogonal basis  of $\R^d$ in which $A$   has
on the diagonal
$0$  or  blocks of the form
$$
\begin{pmatrix}
&0& \; &\lambda_k& \cr
&-\lambda_k& \; &0&
\end{pmatrix},
$$
where $\lambda_k\in\R$ and $i\lambda_k$ is an eigenvalue of $A$.
These are  called  canonical forms.
\item An orthogonal basis  of $\R^d$ for the canonical form can be
constructed in the
following way:
for a fixed eigenvalue $\mu=i\lambda_k\ne 0$ we construct an orthogonal
(orthonormal) basis
$v_1,\dots , v_r$  of the
eigenspace (subspace of $\C^d$) associated to $\mu$, then $\bar
v_1,\dots , \bar v_r$
is  an orthogonal (orthonormal) basis of the
eigenspace associated to $\bar\mu =-\mu$. We put
\begin{equation}\label{canicalbasis}
e_k= \frac{1}{2}(v_k+\bar v_k),\quad f_k= \frac{i}{2}(v_k-\bar v_k).
\end{equation}
Then
$e_1,f_1, \dots, e_r,f_r$ is  an orthogonal (orthonormal)
basis of a real subspace of $\R^d$ of dimension $2r$. In this
basis $A$ has the canonical form.
\end{enumerate}
\end{lemma}

\begin{theorem}\label{diagresanti}
Consider an
analytic family $A:\Omega \subset \R^m \to {\mathcal A}^d$, where
$\mathcal {A}^d$ stands for the space of
antisymmetric $d\times d$ matrices with real entries.
    Then, there
exists
    $ \sigma: W \to
\Omega$  a locally finite composition of blowing-ups with smooth (global)
centers, such that for any $w_0\in W$ there is a neighbourhood $U$  such
that  the corresponding family
$ A \circ \sigma |_{U}:U \to
{\mathcal A}^d$ admits a simultaneous
analytic reduction to the canonical form.
That is, there exists an  analytic choice of vectors
$e:U \to (\R^d)^d$
    such that $ e(w)$ is an orthonormal of  basis  of $\R^d$ and
$A(\sigma (w))$ has
on  diagonal $0$ or blocks of the form
$$
\begin{pmatrix}
&0& \; &\lambda_i& \cr
&-\lambda_i& \; &0&
\end{pmatrix},
$$
for all $w \in U$.
\end{theorem}
\vskip 1cm

\begin{proof}
The arguments are essentially the same as in the proof of Theorem
\ref{diagres}. So we will sketch only the main lines of the proof.
First we resolve
the singularities of the discriminant of the characteristic polynomial
$P(x,z),\, x\in \Omega$
of our family. So we may assume that locally the roots of  of
$P(x,z)$ are analytic
functions of $x$. Recall that if $\lambda (w)$ is such a root then
   $-\lambda (w)$ is also a root of  $P(x,z)$.
We construct (as in the solution of the system \eqref{eigenspace3}) an
orthonormal system of vectors
$v_1(w), \dots,v_r(w)$ which depends analytically on $x,$ in such a way that,
for a generic $w$,
$
V(w) = Ker (A\circ\sigma (w)-\lambda(w)\I_d),
$
is generated by $v_1(w),\dots,v_m(w)$. This requires of course
blowing-ups of the
space of parameters. Now by Lemma
\ref{antisym} vectors $\bar v_1(w),\dots,\bar v_m(w)$  form an
orthonormal basis
of the eigenspace associated to $-\lambda(w)$. So applying formula
\eqref{canicalbasis} we obtain locally a canonical basis for $A(w)$
which depends
analytically on $w$.
\end{proof}

%%%%%%%%

\end{document}